\pdfoutput=1
\documentclass[a4paper,reprint,onecolumn,notitlepage,aip,cha,nofootinbib]{revtex4-2}
\usepackage[left=1.5cm,right=1.5cm,top=1cm,bottom=1.5cm,includeheadfoot]{geometry}

\usepackage[english]{babel}
\usepackage[T1]{fontenc}
\usepackage{tgtermes} 
\usepackage{tgheros} 

\usepackage{amssymb}
\usepackage{amsmath}
\usepackage{amsthm}
\usepackage{graphicx}
\usepackage[dvipsnames]{xcolor}
\usepackage{tikz}
\usepackage[shortlabels]{enumitem}
\usepackage[normalem]{ulem}
\PassOptionsToPackage{hyphens}{url} 
\usepackage[breaklinks=true]{hyperref}
\bibpunct{[}{]}{;}{n}{}{} 

\theoremstyle{plain}
\newtheorem{theorem}{Theorem}
\newtheorem{proposition}[theorem]{Proposition}
\newtheorem{lemma}[theorem]{Lemma}
\newtheorem{corollary}[theorem]{Corollary}
\theoremstyle{definition}
\newtheorem{definition}[theorem]{Definition}
\theoremstyle{remark}
\newtheorem{remark}[theorem]{Remark}
\numberwithin{theorem}{section} 

\newcommand{\R}{\mathbb{R}}

\newcommand{\eps}{\varepsilon}
\newcommand{\proxx}{\prox_{\eps f}(x)}

\newcommand \argmin {\mathop\mathrm{arg\,min}}
\newcommand \dom {\mathop\mathrm{dom}}
\newcommand \prox {\mathop\mathrm{prox}}

\begin{document}
\title{
    Adaptation of Moreau--Yosida regularization to the modulus of convexity
}

\author{Markus Penz}
\email{m.penz@inter.at}
\affiliation{Max Planck Institute for the Structure and Dynamics of Matter, Hamburg, Germany}
\affiliation{Department of Computer Science, Oslo Metropolitan University, Oslo, Norway}

\author{Andre Laestadius}
\affiliation{Department of Computer Science, Oslo Metropolitan University, Oslo, Norway}
\affiliation{Hylleraas Centre for Quantum Molecular Sciences, Department of Chemistry, University of Oslo, Oslo, Norway}

\begin{abstract}
    We study a generalization of Moreau--Yosida regularization that is adapted to the geometry of Banach spaces where the dual space is uniformly convex with modulus of convexity of power type. Important properties for regularized convex functions are given, in particular strong monotonicity of the subdifferential of their convex conjugate and Hölder-continuity of their gradient.
\end{abstract}
\maketitle

\section{Introduction} \label{sec:intro}

Consider the generalization of the Moreau--Yosida (MY) regularization of a convex function on a real Banach space $X$ (see also Definition~\ref{def:MYreg})
\begin{equation}
	f_\eps(x) = \inf \left\{ f(y) + \frac{1}{p\eps}\|x-y\|^p \;\middle|\; y\in X \right\}
\end{equation}
that we call the $p$MY regularization of $f$. The case $p=2$ corresponds to the standard definition of MY regularization. It was demonstrated recently by \citet{Bacho2023} that many of the usual properties of MY regularization carry over to $p$MY regularization if one shifts to a corresponding generalization of the duality map. If $X$ is uniformly smooth with the modulus of smoothness of power type (Definition~\ref{def:p-convex}; this is equivalent to $X^*$ uniformly convex of power type), a useful result of \citet{Xu1991} is applicable that leads to a monotonicity property that will be called ``strong $p^*$-monotonicity'' in Definition~\ref{def:strongly-p-monotone}.

Our interest in this topic comes from the mathematical analysis of density-functional theory~\cite{eschrig2003-book}, a highly popular approximation method in quantum chemistry. This method is typically implemented by using an iteration scheme, where convergence is so far only secured for a MY-regularized version and in finite dimensional spaces~\cite{penz2019guaranteed,penz2020erratum}. The infinite-dimensional generalization leads to a Banach-space formulation using $L^p$ spaces~\cite{KSpaper2018,CDFT-paper} and will possibly rely on the strong-monotonicity result that we prove here.

The structure of the article is as follows. Section~\ref{sec:convex-smooth-spaces} introduces the necessary convexity and smoothness properties of Banach spaces, then Section~\ref{sec:convex-func} repeats basic results from convex analysis, especially the convex conjugate and the subdifferential of a convex function. Section~\ref{sec:duality-map} explains the generalized duality map that will be of special importance in our results, as it is tailored to match the modulus of smoothness. A great summary of many of the topics mentioned before can be found in \citet{schuster-book}.
Then, in Section~\ref{sec:MY} the generalized MY regularization of convex functions with the proximal mapping is introduced and their basic properties are given. Finally, Section~\ref{sec:main-results} lists our main results.
In particular, Theorem~\ref{th:equiv} establishes an equivalence between (i) $p$MY regularization and (ii) adding $(\eps^{p^*-1}/p^*)\|\cdot\|^{p^*}$ to the conjugate function, and this in turn implies (iii) strong $p^*$-monotonicity of the subdifferential of the convex conjugate of the regularized function. As a corollary, Hölder-continuity of the gradient of a regularized function follows.

\section{Convexity and Smoothness of Spaces}
\label{sec:convex-smooth-spaces}

For the convenience of the reader we will repeat several results that are important for the article.
Our basic setting is that of a real Banach space $X$. The topological dual is denoted as $X^*$ and for $x\in X, x^*\in X^*$ the dual pairing $\langle x^*,x \rangle$ simply means applying $x^*$ on $x$, $\langle x^*,x \rangle = x^*(x)$. The norm of these spaces will be denoted as $\|\cdot\|$ in both cases. We recall that $X$ is strictly convex (sometimes also called `rotund') if the unit ball in $X$ is strictly convex, i.e., for any $x,y\in X$, $x\neq y$, $\|x\|=\|y\|=1$, the middle point is strictly inside the ball,
\begin{equation}
    \frac{1}{2}\left\| x+y \right\| < 1.
\end{equation}
Intuitively, strict convexity of the unit ball links to uniqueness of its tangents, and one can make this relation precise on the dual side. A Banach space $X$ is called smooth if for every $x\in X$, $\|x\|=1$, there is exactly one $x^*\in X^*$ such that $\langle x^*,x \rangle=\|x^*\|=1$. This already closely links to the definition of the duality map that is given in a generalized form in Definition~\ref{def:dualitymap} below. Now it holds that $X$ is strictly convex (smooth) if $X^*$ is smooth (strictly convex)~\cite[Theorem~1.101]{Barbu-Precupanu}.
As a refinement we have uniform convexity of $X$ if for any $\eps\in (0,2]$ there is a $\delta(\eps)>0$ such that for all $x,y\in X$ with $\|x\|=\|y\|=1$ and $\|x-y\|=\eps$ we have $1-\frac{1}{2}\|x+y\|\geq \delta(\eps)$. This motivates the following definitions.

\begin{definition}
The \emph{modulus of convexity} of a normed linear space $X$ is defined as the function $\delta_X : (0,2] \to [0,1]$,
\begin{equation}
    \delta_X(\eps) = \inf\left\{ 1-\frac{1}{2}\left\| x+y \right\| \;\middle|\; x,y\in X, \|x\|=\|y\|=1, \|x-y\|=\eps \right\}.
\end{equation}
On the dual side the \emph{modulus of smoothness} of a normed linear space $X$ is the function $\rho_X : [0,\infty) \to [0,\infty)$,
\begin{equation}
    \rho_X(\tau) = \sup\left\{ \frac{1}{2} \left( \| x+y \| + \| x-y \| \right) - 1 \;\middle|\; x,y\in X, \|x\|=1, \|y\|=\tau \right\}.
\end{equation}
\end{definition}

\begin{remark}
In \citet[Sec.~1.4 and 2.2]{chidume-book} the equivalence of various different definitions of $\delta_X$ and $\rho_X$ that can be found in the literature is demonstrated. The characterizations of uniform convex and uniform smoothness are then: $X$ is uniformly convex if and only if $\delta_X(\eps) > 0$ for all $\eps \in (0,2]$. $X$ is uniformly smooth if and only if $\lim_{\tau \to 0}\rho_X(\tau)/\tau = 0$ (this property usually serves as the definition of uniform smoothness).
A Banach space $X$ is uniformly convex (smooth) if and only if $X^*$ uniformly smooth (convex)~\cite[Theorem~1.114]{Barbu-Precupanu}. By the Milman--Pettis theorem, uniform convexity of $X$ (or $X^*$) implies reflexivity of the space. Examples of uniformly convex spaces are the Lebesgue spaces $L^p$, $1<p<\infty$, yet not $L^1, L^\infty$ that are not even strictly convex, nor reflexive. 
The given characterization of uniform convexity (smoothness) can be further narrowed in order to give a definition of even `more' convex (smooth) spaces.
\end{remark}

\begin{definition}\label{def:p-convex}
Let the parameters $p,q > 1$ be real numbers.
We call $X$ a \emph{$p$-convex space} (also: \emph{$p$-uniformly convex}; \emph{uniformly convex (with the modulus of convexity) of power type $p$}) if it has $\delta_X(\eps) \geq c\eps^p$ for some $c>0$. Respectively, $X$ is called a \emph{$q$-smooth space} (also: \emph{$q$-uniformly smooth}) if it has $\rho_X(\tau) \leq c\tau^q$ for some $c>0$.
\end{definition}

\begin{remark}
The terminology already makes it clear that any $p$-convex space is also uniformly convex, which in turn implies strict convexity. This holds analogously for the smoothness properties of spaces. Moreover, it is easily checked that a $p$-convex ($q$-smooth) space is also $p'$-convex ($q'$-smooth) for all $1<p<p'<\infty$ ($1<q'<q<\infty$). In simple terms, for $p$-convexity a lower value of $p$ is the stronger statement, while for $q$-smoothness it is a higher one.
\end{remark}

From this point on, $p$ will fulfill, if not otherwise stated, $1<p<\infty$. Let $p^*$ always denote the Hölder conjugate of $p$, i.e., the number $1<p^*<\infty$ that obeys $1/p+1/p^*=1$.

\begin{proposition}[\citet{chidume-book}, Proposition~5.6]
A space $X$ is $p$-convex ($q$-smooth) if and only if $X^*$ is $p^*$-smooth ($q^*$-convex).
\end{proposition}

That such $p$-convex ($q$-smooth) spaces are by no means rare is demonstrated by the following examples.

\begin{remark}[\citet{chidume-book}, Example~4.13, and \citet{prus1987strongly}]\label{rem:unif-convex}
For $1<p<\infty$ the sequence spaces $\ell^p$, the Lebesgue spaces $L^p$, and the Sobolev spaces $W^{k,p}$ for a domain $\Omega \subseteq \R^n$ are all $\max\{2,p\}$-convex \emph{and} $\min\{2,p\}$-smooth. Hilbert spaces are 2-convex and 2-smooth.
\end{remark}

The remark above indicates that the parameter in the definition of $p$-convexity ($q$-smoothness) can be restricted to $p\geq 2$ ($q\leq 2$) and some references~\cite{pisier1975martingales,prus1987strongly} do this right from the start. To make this restriction obvious, we prove the following lemma from the fact that the modulus of convexity of any normed linear space $X$ is dominated by the modulus of convexity of a Hilbert space (or abstract Euclidean space)~\cite{Nordlander1960}.
The modulus of convexity is the same for any Hilbert space and follows easily from the parallelogram law as $\delta_H(\eps)=1-\sqrt{1-\eps^2/4}$.

\begin{lemma}\label{lem:p2}
Any $p$-convex space always has $p\geq 2$, while a $q$-smooth space always has $q\leq 2$.
\end{lemma}

\begin{proof}
Let $X$ be a $p$-convex space, thus $\delta_X(\eps) \geq c\eps^p$ for some $c>0$.
From \citet{Nordlander1960} we have that generally
\begin{equation}
    \delta_X(\eps) \leq 1-\sqrt{1-\frac{\eps^2}{4}}.
\end{equation}
Both inequalities combine into ($\eps>0$)
\begin{equation}
    \eps^{p-2}(2c-c^2\eps^p) \leq \frac{1}{4}
\end{equation}
and we see that as $\eps\to 0$ this inequality can only be fulfilled if $p\geq 2$. The statement for $q$-smooth spaces then follows from duality.
\end{proof}

\section{Convex functions}
\label{sec:convex-func}

\begin{definition}\label{def:closed-convex}
A convex function $f : X\to \mathbb R\cup \{+\infty\}$ is said to be proper if not identical to $+\infty$, i.e., it has a non-empty domain $\dom{f}:= \{x\in X \mid f(x)\neq +\infty \}\neq \emptyset$. Let $\Gamma_0(X)$ denote the set of proper, convex, lower semi-continuous functions $X \to \mathbb R\cup \{+\infty\}$. On the dual side, $\Gamma_0^*(X^*)$ is the set of proper, convex, weak-* lower semi-continuous functions $X^* \to \mathbb R\cup \{+\infty\}$.
\end{definition}

\begin{lemma}\label{lem:reflexive-gamma}
If $X$ is reflexive, then $\Gamma_0^*(X^*) = \Gamma_0(X^*)$.
\end{lemma}

\begin{proof}
Weak-* semi-continuity is equivalent to weak semi-continuity for $X$ reflexive. But weak \mbox{(semi-)}continuity of convex functions always implies strong (semi-)continuity, so indeed $\Gamma_0^*(X^*) = \Gamma_0(X^*)$.
\end{proof}

\begin{definition}
The \emph{subdifferential} of a convex function $f: X \to \mathbb R \cup \{+ \infty \}$ is defined as
\begin{align}
	\partial f(x) = \{ x^* \in X^* \mid 
	\forall y\in X: f(x)-f(y) \leq \langle x^*,x-y\rangle \}.\label{eq:sub-diff}
\end{align}
\end{definition}

\begin{definition}
A mapping $T : X \to X^*$ is called \emph{monotone} if
\begin{align}
\langle T(x) - T(y), x-y \rangle \geq 0
\end{align}
holds for all $x,y \in X$. If $x \neq y$ and the inequality is strict, the property is accordingly called \emph{strict monotonicity}. The same definition is applicable to set-valued mappings if $T(x)$ is replaced by all $x^* \in T(x)$ when $T(x)$ is non-empty and equivalently for $T(y)$.
\end{definition}

\begin{remark}\label{rem:max-monotone}
It is well known that the subdifferential is always monotone. For $f \in \Gamma_0(X)$ it is even \emph{maximal monotone} \cite{rockafellar1970}, which means it cannot be extended any further in a monotonous fashion. Strict convexity of a function yields strict monotonicity of its subdifferential and we will introduce a special, even stronger form of monotonicity in Definition~\ref{def:strongly-p-monotone} below. Conversely, for a (set-valued) mapping to be equal to the subdifferential of a convex function, the stronger property of \emph{cyclical monotonicity} is needed~\cite[Theorem~2.46]{Barbu-Precupanu}.
\end{remark}

\begin{definition}\label{def:conjugates}
For any function $f: X \to \mathbb R \cup \{\pm \infty \}$ we define the convex conjugate as 
\begin{align}
f^*(x^*) &= \sup \{ \langle x^*,x\rangle - f(x) \mid x \in X \}. \label{eq:def-conjugate}
\end{align}
\end{definition}

\begin{remark}\label{rem:conj-bijection}
With respect to the convex conjugate the sets $\Gamma_0(X)$ and $\Gamma^*_0(X^*)$ have a special relevance: Firstly, application of the convex conjugate on any function always yields functionals of type $\Gamma^*_0(X^*)$, and secondly, the convex conjugate acts as a bijection on those sets with the inverse map just given by the convex conjugate again~\cite[Theorem~2.22]{Barbu-Precupanu}.
\end{remark}


We further note an extremely useful relation between the existence of an optimizer to \eqref{eq:def-conjugate} and the subdifferentials of $f$ and $f^*$.

\begin{theorem}\label{th:inverse}
    Let $f \in \Gamma_0(X)$, then the following are equivalent \cite[Proposition~2.33]{Barbu-Precupanu}:
    \begin{enumerate}[(i)]
        \item $x^* \in \partial f(x)$
        \item\label{th:inverse:superdiff} $x \in \partial f^*(x^*)$
        \item $f(x) + f^*(x^*) = \langle x^*,x \rangle$
    \end{enumerate}
\end{theorem}

\begin{remark}\label{rmk:inverse}
Note that a direct consequence of (i) and (ii) from the previous theorem is that the subdifferential of $f$ is the inverse of the subdifferential of the conjugate of $f$, i.e., $\partial f = [\partial f^*]^{-1}$.
\end{remark}

\section{Generalized Duality Map}
\label{sec:duality-map}

\begin{definition} \label{def:dualitymap}
	The generalized duality map $J_p:X\to X^*$ is
	\begin{equation}\label{eq:dualitymap}
	J_p(x) = \{ x^* \in X^* \mid \|x^*\| = \|x\|^{p-1}, \langle x^*,x \rangle = \|x\|^p \}.
	\end{equation}
\end{definition}

\begin{remark}\label{rem:duality}
If $p=2$ then $J_p=J_2=J$ is the normalized duality map. The generalized duality map is just the subdifferential of the function $\phi_p = \|\cdot\|^p/p$ on $X$, $\partial \phi_p = J_p$ \cite[Proposition~4.9]{chidume-book}. Since further $J_p(x) = \|x\|^{p-2} J(x)$ for $x\neq 0$, which follows directly from the definition, many of the properties of $J$ translate directly to $J_p$. Important properties of the duality map are collected in \citet[Prop.~1.117]{Barbu-Precupanu}. In the present context we need to know that for $X$ reflexive, strictly convex \emph{and} smooth, the duality map $J_p$ is single-valued, demicontinuous (norm-to-weak* continuous), bijective, and has a single-valued inverse that is given by $J_p^{-1} = J_{p^*} : X^* \to X$ \cite{aibinu-mewomo2021}.
Since the same notation $J_{q}$ for a map $X\to X^*$ and a map $X^*\to X$ could cause confusion, the duality map $X^*\to X$ will always be denoted $J_{p}^{-1}$.
The generalized duality map $J_p$ can be better suited than $J$ to reflect the geometry of the underlying Banach space, as it is the case with $L^p$. For the case $X=L^p(\Omega)$ over a measure space $\Omega$, \citet{penot2001yosida} note that it has the simple form $J_p(x)(\omega) = |x(\omega)|^{p-2} x(\omega)$ for almost all $\omega \in \Omega$ that does not involve any integration as it would be the case for the normalized duality map, which can be checked by direct comparison with the conditions in \eqref{eq:dualitymap}. Note that in a Hilbert space setting, where the duality map $J$ is precisely the canonical isomorphism given by the Riesz representation theorem, $J$ must in general \emph{not} be the identity function, as shown with the example of Sobolev spaces $H^1_0$ and $H^{-1}$~\cite[Section~5.9.1]{evans-book}. However, $J$ is always just the identity if $X$ is identified with its dual $X^*$, like it will be the case for $X = \ell^2$ or $L^2$.
\end{remark}

We can connect the $p$-convexity of a space to a feature of its generalized duality map that we suggest to call strong $p$-monotonicity.

\begin{definition}\label{def:strongly-p-monotone}
A mapping $T : X \to X^*$ is called \emph{strongly $p$-monotone} if there is a $c>0$ such that for all $x,y \in X$ it holds
\begin{align}
\langle T(x) - T(y), x-y \rangle \geq c \|x-y\|^{p}.
\end{align}
\end{definition}

\begin{lemma}[\citet{Xu1991}, Corollary 1 (i) and (iii); also \citet{chidume-book}, Corollary 4.17 (i) and (iii)]\label{lem:Xu}
    Let $p\geq 2$, then $X$ is $p$-convex if and only if the generalized duality map $J_p(y)$ is strongly $p$-monotone.
\end{lemma}

\begin{lemma}\label{lem:duality-Hoelder}
Let $X$ be strictly convex and $p$-smooth, $p\leq 2$, then $J_p$ is $(p-1)$-Hölder continuous. If $p=2$ then this implies Lipschitz continuity.
\end{lemma}

\begin{proof}
That $X$ is $p$-smooth implies reflexivity and $p^*$-convexity for $X^*$. The conditions imply that $J_p^{-1}=J_{p^*}$ exists and is single-valued.
Now, Lemma~\ref{lem:Xu}, but with the roles of $X$ and $X^*$ interchanged, shows that there is a $c > 0$ such that for all $x^*,y^*\in X^*$ (strong $p^*$-monotonicity of $J_p^{-1}=J_{p^*}$)
\begin{equation}
    \langle x^*-y^*,J_p^{-1}(x^*)-J_p^{-1}(y^*) \rangle \geq c \|x^*-y^*\|^{p^*}.
\end{equation}
Here, the strict convexity of $X$ guarantees that $J_p^{-1}$ is single-valued.
Inserting $x=J_p^{-1}(x^*), y=J_p^{-1}(y^*)$ and estimating the left-hand side above yields
\begin{equation}
    c \|J_p(x)-J_p(y)\|^{p^*} \leq \langle J_p(x)-J_p(y),x-y \rangle \leq \|J_p(x)-J_p(y)\| \|x-y\|
\end{equation}
and we thus have
\begin{equation}
    \|J_p(x)-J_p(y)\| \leq \left(\frac{1}{c}\|x-y\|\right)^{1/(p^*-1)} = \left(\frac{1}{c}\|x-y\|\right)^{p-1}.
\end{equation}
This shows that the generalized duality map is $(p-1)$-Hölder continuous.
\end{proof}

\begin{remark}
This result explains the relevance of Lemma~\ref{lem:p2} that showed that a $p$-convex space always has $p\geq 2$. In the context of the lemma above it is thus $p^*\geq 2$, which implies that the Hölder conjugate has $p\leq 2$ and thus $p-1\leq 1$ as it must be for Hölder continuity. This is because an $\alpha$-Hölder continuous function $F:X\to Y$ with $\alpha>1$ can only be constant. This can be seen by taking $x,y\in X$ with $x\neq y$, then the Hölder property means that
\begin{equation}
    \frac{\|F(x)-F(y)\|}{\|x-y\|}\leq C \|x-y\|^{\alpha-1}.
\end{equation}
So the difference quotient for $F$ approaches zero as $x\to y$. Hence the Fréchet derivative of $F$ exists and is zero everywhere, making $F$ a constant function.
\end{remark}

\section{Generalized Moreau--Yosida regularization}
\label{sec:MY}

We open this section with a basic result about strict convexity. Remember that a function $f:X\to \R\cup\{+\infty\}$ is called strictly convex if for any $x,y\in X$, $x\neq y$, and all $\lambda\in(0,1)$ it holds
\begin{equation}\label{eq:func-strict-convex}
    f(\lambda x+(1-\lambda)y) < \lambda f(x) + (1-\lambda)f(y).
\end{equation}
The strict convexity of $\|x\|^p$, $p>1$, relates to strict convexity of the space itself.

\begin{lemma}\label{lem:PhiStrictConvex}
    For all $p>1$, the function $\phi_p : X \to \R, x \mapsto \|x\|^p/p$ is strictly convex if and only if the space $X$ is strictly convex.
\end{lemma}

\begin{proof}
Assume $X$ strictly convex. Take any $x,y \in X$, $x\neq y$, $\lambda \in (0,1)$, then by the triangle inequality $\|\lambda x + (1-\lambda)y\|^p \leq (\lambda \|x\| + (1-\lambda)\|y\|)^p$. Since $t \mapsto t^p$ is clearly strictly convex on $\R$ for $p>1$, it follows $\|\lambda x + (1-\lambda)y\|^p < \lambda \|x\|^p + (1-\lambda) \|y\|^p$ if $\|x\|\neq \|y\|$. In the case $\|x\| = \|y\|$, set $\|x\| = \|y\| = 1$ without loss of generality (since we can always rescale $x/\|x\|$ and $y/\|y\|$ and use $\|x\| = \|y\|$) and then resort to strict convexity of $X$ to show $\|\lambda x + (1-\lambda)y\|^p < \lambda \|x\|^p + (1-\lambda) \|y\|^p$ for $x\neq y$ from $\|\lambda x + (1-\lambda)y\| < 1$.\\
Conversely, the strict convexity of $X$ follows directly from the defining inequality of strict convexity of $\phi_p$ when setting $\|x\| = \|y\| = 1$.
\end{proof}

A function $f:X\to \R\cup\{+\infty\}$ is called uniformly convex, if in the condition for convexity an additional modulus $\varphi:\R_{\geq 0} \to \R_{\geq 0}$ that has $\varphi(\tau)=0$ only at $\tau=0$ can be subtracted,
\begin{equation}
    f(tx+(1-t)y) \leq tf(x) + (1-t)f(y)-t(1-t)\varphi(\|x-y\|).
\end{equation}
A similar result as above can then be stated for $p$-convexity of spaces.

\begin{lemma}[\citet{chidume-book}, Theorem~4.16]
    For all $p>1$, the function $\phi_p : X \to \R, x \mapsto \|x\|^p/p$ is uniformly convex if and only if the space $X$ is $p$-convex.
\end{lemma}

\begin{lemma}\label{lem:PhiGateauxFrechet}
	If $X$ is smooth then $\phi_p$ is Gâteaux-differentiable. If $X$ is uniformly smooth then $\phi_p$ is Fréchet-differentiable. In both cases, the derivative of $\phi_p$ is the generalized duality map $J_p$.
\end{lemma}

\begin{proof}
That the subdifferential of $\phi_p$ is identical to the generalized duality map $J_p$ was already mentioned in Remark~\ref{rem:duality}. If $X^*$ is strictly convex then the normalized duality map is single-valued \citep[Proposition 1.117(iii)]{Barbu-Precupanu}, which because of $J_p(x) = \|x\|^{p-2} J(x)$ also means that $J_p$ is single-valued. This in turn already establishes Gâteaux differentiability \citep[Proposition 2.40]{Barbu-Precupanu}. If furthermore $X^*$ is uniformly convex then the normalized duality map is also continuous \citep[Proposition~1.117(vi)]{Barbu-Precupanu}.
Continuity of $J_p$ is not lost close to the origin, as can be seen using $J$ homogeneous by writing $J_p(x) = \|x\|^{p-2} J(\|x\|\cdot x/\|x\|) = \|x\|^{p-1} J(x/\|x\|)$.
Finally, continuity of the differential implies Fréchet differentiability \citep[Lemma~34.3]{Blanchard-Bruening}.
\end{proof}

\begin{definition}\label{def:MYreg}
    For $p>1$ and $\eps>0$, the $p$-Moreau--Yosida ($p$MY) regularization of $f \in \Gamma_0(X)$ is defined as the \emph{infimal convolution} of $f$ with $\eps^{-1}\phi_p$,
	\begin{equation}\label{eq:MYreg}
	f_\eps(x) = \inf \{ f(y) + \eps^{-1}\phi_p(x-y) \mid y\in X \}.
	\end{equation}
\end{definition}

\begin{figure}[ht]
\centering
\begin{tikzpicture}[scale=.65]
    \draw[->] (-1,0) -- (10,0)
        node[right] {$x$};
    
    \fill (0,0) circle [radius=2pt] node[below=5pt] {$x_0$}; 
    \fill (1,0) circle [radius=2pt] node[below=5pt] {$x_1$}; 
    \fill (2,0) circle [radius=2pt] node[below=5pt] {$x_2$}; 
    \fill (3,0) circle [radius=2pt] node[below=5pt] {$x_3$}; 
    \fill (4,0) circle [radius=2pt] node[below=5pt] {$x_4$}; 
    
    \draw[black!30] (-1,6) parabola bend (1,4) (3,6);
    \draw[black!30] (-.45,6) parabola bend (2,3) (4.45,6);
    \draw[black!30] (.17,6) parabola bend (3,2) (5.83,6);
    \draw[black!30] (.84,6) parabola bend (4,1) (7.16,6);
    
    \draw (-1,5.5) -- (3.5,1);
    \draw (3.5,1) parabola bend (4.5,.5) (5.5,1);
    \draw (5.5,1) -- (10,5.5) node[right] {$f_\eps$};
    
    \draw[dashed] (-1,6) -- (4.5,.5);
    \draw[dashed] (4.5,.5) -- (10,6) node[right] {$f$};
\end{tikzpicture}
\caption{Moreau--Yosida regularization $f_\eps$ (solid) with $p=2$ of an exemplary, one-dimensional $f$ (dashed), showing also the regularization parabolas $\frac{1}{2\eps}\|x-x_i\|^2+\text{const}$, $i\geq 1$, that trace out $f_\eps$ (grey). The sequence $\{ x_i \}_i$ approaches the minimum following the proximal-point algorithm.}
\label{fig:MY}
\end{figure}
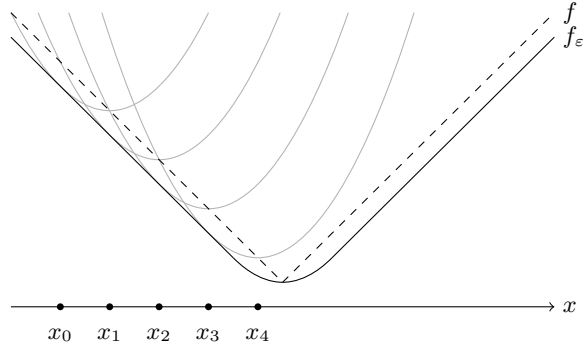

\begin{remark}\label{rem:MY1}
The case $p=2$ recovers the standard MY regularization. For reasons of brevity, we will not include $p$ in the notation for the $p$MY regularization. This generalized form of the MY regularization can already be found in the textbook \citet[Proposition~12.15]{Bauschke-Combettes}, where it is defined on Hilbert spaces. \citet{Bacho2023} recently showed that the main properties of the standard MY regularization also hold in the case of $p$MY and gave a great literature review. \citet{kuwae2015resolvent} and \citet{bacak-kohlenbach-2018} present this generalized MY regularization with $\eps^{-1}$ replaced by $\eps^{1-p}$ under the name \emph{Hopf--Lax formula} that is related to solutions of a Hamilton--Jacobi equation and consequently $f_\eps$ is also called a \emph{Hamilton--Jacobi semigroup} (in $t=\eps^{1/(p-1)}$).
\citet{penot2001yosida} and \citet{bacak-kohlenbach-2018} allow for the even more general \emph{Young functions} instead of just $\phi_p$.
Further, in \citet{penot2001yosida} a Yosida approximation for set-valued mappings from $X$ into $X^*$ is discussed alongside its relation to the MY regularization of functionals. Another form of generalization of the standard MY regularization, by replacing $\phi_p(x)$ with $\|Ax\|^2/2$, is briefly discussed in the conclusions of \citet{Kvaal2014}. Therein, some interesting, possible connections to approaches employed in density-functional theory that resemble the regularization technique are noted as well.
The definition of regularization as an infimal convolution has the benefit of yielding convenient relations between the convex conjugates of regularized and unregularized functionals, expressed through equations such as \eqref{eq:feconj} and \eqref{eq:feconj-superdiff}.
\end{remark}

\begin{remark}\label{rem:MY2}
\citet{penot2001yosida} show that (as expected) $f_\eps \in \Gamma_0(X)$, more specifically that the regularized $f_\eps$ is even a continuous, convex function \cite[Proposition~5.6]{penot2001yosida}. From the definition it follows directly that $f_\eps \leq f$ everywhere and that for all $x,y \in X$ it holds
$f_\eps(x) \leq f(y) + \eps^{-1}\phi_p(x-y)$, 
which in the context of Fig.~\ref{fig:MY} means that the regularization paraboloids are always above $f_\eps$.
Further, for any minimum $x\in X$ of $f$ it holds that also $f_\eps$ has a minimum at the same $x$. The infimum in \eqref{eq:MYreg} above is always uniquely attained, because the added $\phi_p$ is strictly convex if $X$ is strictly convex by Lemma~\ref{lem:PhiStrictConvex} \cite[Proposition~1.103(iv), Theorem~2.11, Remarks~2.12 and 2.13]{Barbu-Precupanu}. This unique minimizer for any given $x$ allows for the definition of the \emph{proximal mapping}.
\end{remark}

\begin{definition}[proximal mapping]\label{def:prox}
    Let $X$ be strictly convex, $f \in \Gamma_0(X)$, $p>1$, and $\eps > 0$. For all $x \in X$ the proximal mapping is
    \begin{equation}
    \begin{aligned}
    \proxx &= \argmin\{ f(y) + \eps^{-1}\phi_p(x-y) \mid y\in X \} \\
    &= \argmin\{ \eps f(y) + \phi_p(x-y) \mid y\in X \}.
    \end{aligned}
    \end{equation}
\end{definition}

\begin{remark}
    The above definition leads, by seeking the minimum in \eqref{eq:MYreg}, to the following useful relations,
    \begin{align}
    \label{eq:prox-df}
    &\partial f(\proxx) - \eps^{-1} J_p(x-\proxx) \ni 0, \\
    \label{eq:f-eps-prox}
    &f_\eps(x) = f(\proxx) + \eps^{-1}\phi_p(x-\proxx).
    \end{align}
\end{remark}

\begin{remark}
In the limit $\eps\to 0+$ the proximal mapping approaches the original point. More precisely, let $x\in \dom{f} $ then $\Vert x - \proxx \Vert = \mathcal O(\eps^{1/p})$ as $\eps\to 0+$. We show this statement in the simplified setting of $f\geq 0$ (the general case can be treated by the fact that every $f\in\Gamma_0(X)$ has an affine linear minorant). From \eqref{eq:f-eps-prox} and $f_\eps\leq f$ we have
\begin{equation*}
    f(\proxx) + \eps^{-1}\phi_p(x-\proxx) \leq f(x) =:C.
\end{equation*}
The positivity of $f$ then entails $\Vert x-\proxx \Vert^p \leq Cp\eps $ and the claim follows.
\end{remark}

\begin{remark}[proximal-point algorithm]\label{rem:prox-point-algorithm}
For $f \in \Gamma_0(X)$ a sequence defined by $x_{i+1} = \prox_{\eps f}(x_i)$ gradually approaches the infimum of $f$ as depicted in Fig.~\ref{fig:MY}, see also \citet[Theorem~27.1]{Bauschke-Combettes}. 
Anticipating Theorem~\ref{lemma:frechet}, this amounts exactly to a gradient-descent algorithm for the regularized $f_\eps$.
\end{remark}

The following result taken from \citet{Bacho2023} that we quote only partly gives important properties of the $p$MY regularization, especially the differentiability of regularized functionals. Like with the MY regularization one sees that the differentiability properties from $\phi_p$, noted on in Lemma~\ref{lem:PhiGateauxFrechet}, are inherited. The result for Hilbert spaces is standard~\cite[Corollary~2.59]{Barbu-Precupanu}.

\begin{theorem}\label{lemma:frechet}
	Let $X$ be reflexive, strictly convex and smooth, $p>1$, $\eps>0$, and $f\in \Gamma_0(X)$ with $p$MY regularization $f_\varepsilon$. 
    Then $f_\eps$ is convex and locally Lipschitz-continuous, and if $f$ is strictly convex then $f_\eps$ is as well.
    Further, $f_\eps$ is Gâteaux-differentiable on $X$ with a demicontinuous derivative that satisfies
	\begin{equation}\label{eq:nabla-fe}
	\nabla f_\eps(x)= \eps^{-1}J_p(x - \proxx).
	\end{equation}
	If $X$ is uniformly smooth, then $f_\eps$ is even Fréchet-differentiable. If $X$ is a Hilbert space and $p=2$, then $\nabla f_\eps$ is Lipschitz-continuous with constant $\eps^{-1}$.
\end{theorem}

\begin{remark}
Equation~\eqref{eq:nabla-fe} has an easy pictorial interpretation visible in Fig.~\ref{fig:MY}: Since $J_p$ is just the differential of $\phi_p$, we rewrite $\nabla f_\eps(x)= \eps^{-1}\nabla\phi_p(x - \proxx)$ and see that at any point $x$ the regularized function $f_\eps$ is aligned tangentially to $\eps^{-1}\phi_p$ (centered around the proximal point of $x$).
\end{remark}

\section{Main Results}
\label{sec:main-results}

\begin{theorem}\label{th:equiv}
Let $X$ be strictly convex and smooth, $p>1$ and $\eps>0$, then the following statements are equivalent for $f,g \in \Gamma_0(X)$:
\begin{enumerate}[(i)]
    \item\label{th:equiv:f-eps} 
    $f=g_\eps$ ($f$ is the $p$MY regularization of $g$),
    
    \item\label{th:equiv:f-subtract}
    $f^* = g^* + \eps^{p^*-1}\phi_{p^*}$.
\end{enumerate}
If further $X$ is $p$-smooth ($p\leq 2$) then this implies:
\begin{enumerate}[(i)]
    \setcounter{enumi}{2}
    \item\label{th:equiv:strongly-mon} 
    $\partial f^*$ is strongly $p^*$-monotone by Definition~\ref{def:strongly-p-monotone}, i.e.,
    there is a $c>0$ such that for all $x^*, y^* \in X^*$ and all $x \in \partial f^*(x^*)$, $y \in \partial f^*(y^*)$ with non-empty subdifferentials it holds
    \begin{equation}\label{eq:g-mon}
    \langle x^*-y^*,x-y \rangle \geq c \|x^*-y^*\|^{p^*}.
    \end{equation}
\end{enumerate}
\end{theorem}

\begin{proof}
In the proof we will use, as noted in Remark~\ref{rem:duality}, that under the given assumptions on $X$ the duality map $J_{p}:X \to X^*$ is a bijection with a well-defined and single-valued inverse $J_{p}^{-1}$.

(i) $\Rightarrow$ (ii). We note that by definition
\begin{equation}
\begin{aligned}
	f^*(x^*) = (g_\eps)^*(x^*) =& \sup\{ \langle x^*,x \rangle - \inf\{g(y) + \eps^{-1}\phi_p(x-y) \mid y\in X \} \mid x\in X \} \\
    =& -\inf\{ g(y) + \eps^{-1}\phi_p(x-y) - \langle x^*,x \rangle \mid x,y\in X \} \\
	=& -\inf\{ g(y) - \langle x^*,y \rangle + \eps^{-1}\phi_p(x-y) - \langle x^*,x-y \rangle \mid x,y \in X \} \\
	=& -\inf\{ g(y) - \langle x^*,y \rangle \mid y \in X\} - \inf\{ \eps^{-1}\phi_p(z) - \langle x^*,z \rangle \mid z \in X \} \\
    =& \sup\{ \langle x^*,y \rangle - g(y) \mid y \in X\} + \sup\{ \langle x^*,z \rangle - \eps^{-1}\phi_p(z) \mid z \in X \} \\
	=& g^*(x^*) + (\eps^{-1}\phi_p)^*(x^*).
\end{aligned}
\end{equation}
From $\phi_p(x) =  \| x \|^p/p$ it follows by definition that $\phi_p^*(x^*) = \sup\{ \langle x^*,x \rangle - \| x \|^p/p \mid x\in X \}$, where the maximum is uniquely attained because of strict convexity. Seeking the maximum by subdifferentiation yields $J_p(x)=x^*$, where the duality map is single-valued because of smoothness. The definition of the generalized duality map then provides the relations $\|x^*\|=\|x\|^{p-1}$ and $\langle x^*,x \rangle=\|x\|^p$ that must hold for the maximum. Inserting these two equalities into the expression of the convex conjugate of $\phi_p$ gives $\phi_p^*(x^*) = (1-1/p)\|x\|^p = \|x^*\|^{p^*}/p^* = \phi_{p^*}(x^*)$. Additionally, the scaling relation $(\lambda h)^*(x^*) = \lambda h^*(x^*/\lambda)$, $\lambda>0$, yields $(\eps^{-1}\phi_p)^*(x^*)=\eps^{-1}\|\eps x^*\|^{p^*}/p^*=\eps^{p^*-1} \phi_{p^*}(x^*)$ and thus $f^* = g^* + \eps^{p^*-1}\phi_{p^*}$ which proves (ii).

(ii) $\Rightarrow$ (i). Starting from $f^* = g^* + \eps^{p^*-1}\phi_{p^*}$ we can do all steps above in reverse order and thus conclude that $f^*=(g_\eps)^*$. But this means $f=g_\eps$ by invertibility of the convex conjugate.

(ii) $\Rightarrow$ (iii).
Because of strong $p^*$-convexity of $f^*$ we know there is a $g \in \Gamma_0(X)$ (cf.~Remark~\ref{rem:conj-bijection}) such that
\begin{equation}\label{eq:feconj}
    f^* = g^* + \eps^{p^*-1}\phi_{p^*},
\end{equation}
which entails
\begin{equation}\label{eq:feconj-superdiff}
    \partial f^* = \partial g^* + \eps^{p^*-1} J^{-1}_{p}
\end{equation}
on $X^*$ by $\partial\phi_{p^*} = J_p^{-1}$ (cf.~Remark~\ref{rem:duality}). Note that the sum rule for the subdifferentials only holds because the involved functionals are all convex and $\phi_{p^*}$ is continuous~\cite[Theorem~3.57]{Barbu-Precupanu}.
Now for any $\tilde x \in \partial g^*(x^*)$, $\tilde y \in \partial g^*(y^*)$ by definition of the subdifferential it holds
\begin{align}
    g^*(x^*) - g^*(y^*) &\leq \langle x^*-y^*,\tilde x \rangle,\\
    g^*(y^*) - g^*(x^*) &\leq \langle y^*-x^*,\tilde y \rangle = \langle x^*-y^*,-\tilde y \rangle.
\end{align}
Adding those two inequalities gives (monotonicity of the subdifferential, as already noted in Remark~\ref{rem:max-monotone})
\begin{equation}\label{eq:tilde-mon}
    \langle x^*-y^*,\tilde x-\tilde y \rangle \geq 0.  
\end{equation}
Equation~\eqref{eq:feconj-superdiff} results in the relations $x = \tilde x + \eps^{p^*-1} J_p^{-1}(x^*)$ and $y = \tilde y + \eps^{p^*-1} J_p^{-1}(y^*)$ for $x\in \partial f^*(x^*)$, $y\in \partial f^*(y^*)$, and inserting this into \eqref{eq:tilde-mon} yields
\begin{equation}\label{eq:tilde-mon-2}
    \langle x^*-y^*,x-y \rangle \geq \eps^{p^*-1} \langle x^*-y^*,J_p^{-1}(x^*)-J_p^{-1}(y^*) \rangle.
\end{equation}
From the additional assumption that $X$ is $p$-smooth follows that $X^*$ is $p^*$-convex, so Lemma~\ref{lem:Xu}, again with the roles of $X$ and $X^*$ interchanged, shows that for $X^*$ $p^*$-convex there is a $c_1 > 0$ such that
\begin{equation}
    \langle x^*-y^*,J_p^{-1}(x^*)-J_p^{-1}(y^*) \rangle \geq c_1 \|x^*-y^*\|^{p^*}.
\end{equation}
This together with \eqref{eq:tilde-mon-2} and setting $c=\eps^{p^*-1} c_1$ establishes (iii).
%
\end{proof}

\begin{remark}
Theorem~\ref{th:equiv} has also \ref{th:equiv:strongly-mon} equivalent to \ref{th:equiv:f-eps} and \ref{th:equiv:f-subtract} in the case $X=\R^n$ with $p=p^*=2$, see \citet[Proposition~12.60]{rockafellar2009variational}.
A similar result for Hilbert spaces can be found in \citet[Theorem~18.15]{Bauschke-Combettes}, the strong monotonicity of $\partial f^*$ is then called cocoercivity of $\nabla f$.
When comparing to such results, one is tempted to call the property in \ref{th:equiv:strongly-mon} ``strong $p^*$-convexity'', but note that this notion does not make sense for $p\neq 2$. This is seen easily by taking for example $f(x)=\|x-x_0\|^3$, which should be clearly ``strongly 3-convex'', but subtracting $c\|x\|^3$ with any $c>0$ from $f(x)$ leaves us with a function that is not convex.
\end{remark}


The theorem above now implies Hölder-continuity for the gradient of the $p$MY regularization of an $f\in\Gamma_0(X)$ just like it was shown for the generalized duality map in Lemma~\ref{lem:duality-Hoelder}.

\begin{corollary}
    Let $X$ be strictly convex and $p$-smooth, $p\leq 2$ and $\eps > 0$, then for any $f \in \Gamma_0(X)$ we have that $\nabla f_\eps$ is $(p-1)$-Hölder continuous.
\end{corollary}

\begin{proof}
    From Theorem~\ref{th:equiv}\ref{th:equiv:strongly-mon} we have strong $p^*$-monotonicity of $\partial (f_\eps)^*$ and since $\nabla f_\eps = [ \partial (f_\eps)^* ]^{-1}$ (Theorem~\ref{th:inverse}) we have for $x^*=\nabla f_\eps(x), y^*=\nabla f_\eps(y)$ that
    \begin{equation}
    c \|x^*-y^*\|^{p^*} \leq \langle x^*-y^*,x-y \rangle \leq \|x^*-y^*\| \|x-y\|.
    \end{equation}
    This immediately gives
    \begin{equation}
    \|x^*-y^*\| \leq \left(\frac{1}{c}\|x-y\|\right)^{1/(p^*-1)} =\left(\frac{1}{c}\|x-y\|\right)^{p-1},
    \end{equation}
    which concludes the proof.
\end{proof}

\acknowledgments
MP and AL were supported by the ERC-2021-STG
grant agreement No.~101041487 REGAL. AL was also
supported by the Research Council of Norway through
CoE Hylleraas Centre for Quantum Molecular Sciences
Grant No.~262695 and CCerror Grant No.~287906.

%

\end{document}